%% file: gquiver.tex
\documentclass{amsart}

\usepackage{amsfonts}
\usepackage{amssymb}
\usepackage{amscd}
\usepackage{graphicx}
\usepackage{stmaryrd}
\usepackage{xypic}
\usepackage{psfrag}

\input{preamble}
\input{tableau}
\input{psfrag}

\numberwithin{equation}{section}

\begin{document}

\title{Grothendieck classes of quiver varieties}
\author{Anders Skovsted Buch}
\address{Massachusetts Institute of Technology \\
  Building 2, Room 248 \\
  77 Massachusetts Avenue \\
  Cambridge, MA 02139
}
\date{\today}
\email{abuch@math.mit.edu} 
\thanks{The author was partially supported by NSF Grant DMS-0070479}

\begin{abstract}
  We prove a formula for the structure sheaf of a quiver variety in
  the Grothendieck ring of its embedding variety.  This formula
  generalizes and gives new expressions for Grothendieck polynomials.
  We furthermore conjecture that the coefficients in our formula have
  signs which alternate with degree.  The proof of our formula
  involves $K$-theoretic generalizations of several useful
  cohomological tools, including the Thom-Porteous formula, the
  Jacobi-Trudi formula, and a Gysin formula of Pragacz.
\end{abstract}

\maketitle

\input{intro}
\input{groth}
\input{intseq}
\input{formula}
\input{appgroth}
\input{jacobi-trudi}
\input{gysin}

\input{bibliography}

\end{document}

%% file: preamble.tex
\newcommand{\N}{{\mathbb N}}
\newcommand{\Z}{{\mathbb Z}}

\renewcommand{\P}{{\mathbb P}}

\newcommand{\Groth}{{\mathfrak{G}}}
\newcommand{\bull}{{\sssize \bullet}}
\newcommand{\rank}{\operatorname{rank}}
\newcommand{\sh}{\operatorname{sh}}

\newcommand{\depth}{\operatorname{depth}}
\newcommand{\Tor}{\operatorname{Tor}}
\newcommand{\Hom}{\operatorname{Hom}}
\newcommand{\Gr}{\operatorname{Gr}}
\newcommand{\Fl}{\operatorname{F\ell}}
\renewcommand{\O}{{\mathcal O}}
\newcommand{\OOmega}{\boldsymbol{\Omega}}

\newcommand{\bbn}[2]{\left[\begin{smallmatrix} #1 \\ #2 \end{smallmatrix}\right]}

\newcommand{\pic}[2]{\includegraphics[scale=0.#1]{Figures/#2.eps}}

\newcommand{\picB}[1]{\includegraphics[scale=0.90]{Figures/#1.eps}}
\newcommand{\picC}[1]{\includegraphics[scale=0.60]{Figures/#1.eps}}
\newcommand{\picD}[1]{\includegraphics[scale=0.50]{Figures/#1.eps}}

\newtheorem{thm}{Theorem}[section] 
\newtheorem{lemma}[thm]{Lemma} 

\newtheorem{cor}[thm]{Corollary}
\newtheorem{conj}[thm]{Conjecture}

\theoremstyle{definition}
\newtheorem{example}[thm]{Example}
\newtheorem{defn}[thm]{Definition}

\newcommand{\refsec}[1]{Section~\ref{#1}}

\newcommand{\refthm}[1]{Theorem~\ref{#1}}

\newcommand{\reflemma}[1]{Lemma~\ref{#1}}

\newcommand{\refconj}[1]{Conjecture~\ref{#1}}

\newcommand{\refexm}[1]{Example~\ref{#1}}

\newenvironment{romenum}{\begin{enumerate}}{\end{enumerate}}

\newcommand{\lab}[1]{\label{#1}}


%% file: tableau.tex
\newlength{\cellsize} 
\newsavebox{\cell}
\sbox{\cell}{\begin{picture}(18,18)
\put(0,0){\line(1,0){18}}
\put(0,0){\line(0,1){18}}
\put(18,0){\line(0,1){18}}
\put(0,18){\line(1,0){18}}
\end{picture}}
\newcommand\cellify[1]{\def\thearg{#1}\def\nothing{}%
\ifx\thearg\nothing
\vrule width0pt height\cellsize depth0pt\else
\hbox to 0pt{\usebox{\cell} \hss}\fi%
\vbox to \cellsize{
\vss
\hbox to \cellsize{\hss$#1$\hss}
\vss}}
\newcommand\tableau[1]{\vtop{\let\\\cr
\baselineskip -16000pt \lineskiplimit 16000pt \lineskip 0pt
\ialign{&\cellify{##}\cr#1\crcr}}}


%% file: psfrag.tex
\psfrag{ay1}{$\mu$}
\psfrag{ay2}{$\lambda$}

\psfrag{rij-1-rij}{$r_{i,j-1} - r_{ij}$}
\psfrag{ri+1j-rij}{$r_{i+1,j} - r_{ij}$}
\psfrag{R01}{$R_{01}$}
\psfrag{R12}{$R_{12}$}
\psfrag{R23}{$R_{23}$}
\psfrag{Rn-1n}{$R_{n-1,n}$}
\psfrag{R02}{$R_{02}$}
\psfrag{R13}{$R_{13}$}
\psfrag{Rn-2n}{$R_{n-2,n}$}
\psfrag{sig1}{$\sigma_1$}
\psfrag{sig2}{$\sigma_2$}
\psfrag{tau1}{$\tau_1$}
\psfrag{taun-1}{$\tau_{n-1}$}
\psfrag{Ri-1i}{$R_{i-1,i}$}
\psfrag{sigi}{$\sigma_i$}
\psfrag{taui-1}{$\tau_{i-1}$}


%% file: intro.tex
\section{Introduction}

Let $X$ be a non-singular variety and $E_0 \to E_1 \to \dots \to
E_n$ a sequence of vector bundles and bundle maps over $X$.  A set of
{\em rank conditions\/} for this sequence is a collection $r =
\{r_{ij}\}$ of non-negative integers for $0 \leq i < j \leq n$.  This
data defines a {\em quiver variety\/}:
\begin{equation} \label{eqn:locus}
  \Omega_r = \Omega_r(E_\bull) = \{ x \in X \mid 
   \rank(E_i(x) \to E_j(x)) \leq r_{ij} ~\forall i < j \} \,.
\end{equation}
This set has a natural structure of subscheme of $X$.  Namely, it is
the scheme theoretic intersection of the zero sections of the bundle
maps $\bigwedge^{r_{ij}+1} E_i \to \bigwedge^{r_{ij}+1} E_j$. 

We will demand that the rank conditions can {\em occur}, i.e.\ there
exists a sequence of vector spaces and linear maps $V_0 \to V_1 \to
\cdots \to V_n$ such that $\dim V_i = \rank E_i$ and $\rank(V_i \to
V_j) = r_{ij}$ for all $i < j$.  If we set $r_{ii} = \rank E_i$, then
this is equivalent to the conditions $r_{ij} \leq \min(r_{i,j-1},
r_{i+1,j})$ for $i < j$ and $r_{i+1,j-1} - r_{i,j-1} - r_{i+1,j} +
r_{ij} \geq 0$ for $j - i \geq 2$.

The expected (and maximal possible) codimension of the quiver variety
$\Omega_r$ is $d(r) = \sum_{i < j} (r_{i,j-1} - r_{ij}) (r_{i+1,j} -
r_{ij})$.  When this codimension is obtained, the main result of
\cite{buch.fulton:chern} gives a formula for the cohomology class of
$\Omega_r$:
\begin{equation} \label{eqn:cohom_quiver}
   [\Omega_r] = \sum_{|\mu| = d(r)} c_\mu(r) \, s_{\mu_1}(E_1 - E_0) \,
   s_{\mu_2}(E_2 - E_1) \cdots s_{\mu_n}(E_n - E_{n-1}) \,.
\end{equation}
This sum is over sequences $\mu = (\mu_1, \dots, \mu_n)$ of $n$
partitions such that the sum $|\mu| = \sum |\mu_i|$ of the weights of
these partitions is equal to the expected codimension $d(r)$.  (Recall
that the weight of a partition is the sum of its part, or the number
of boxes in its Young diagram.)  If $\lambda$ is a partition then
$s_\lambda(E_i - E_{i-1})$ denotes the double Schur function
$s_\lambda(x;y)$ applied to the Chern roots of the bundles $E_i$ and
$E_{i-1}$.  The coefficients $c_\mu(r)$ are certain integers given by
an explicit combinatorial algorithm.  Surprisingly, these coefficients
appear to be non-negative.  It is conjectured in
\cite{buch.fulton:chern} that each coefficient counts the number of
sequences of semistandard Young tableaux satisfying certain
properties; this has been proved when the sequence $E_\bull$ has at
most four bundles \cite{buch:on}.

While the cohomology class of a quiver variety $\Omega_r$ represents
useful global information, there is more information hidden in its
structure sheaf $\O_{\Omega_r}$.  The best possible representation of
this information that one could hope for might be an explicit
resolution of the structure sheaf by locally free sheaves on $X$.
Such a resolution would generalize fundamental constructions such as
the Koszul complex and the Eagon-Northcott complex
\cite{eagon.northcott:ideals}, at least up to quasi-isomorphism.  Such
resolutions, however, are known in only very few cases, such as for
Schubert varieties in Grassmannians \cite{lascoux:foncteurs}.  The
main theorem in this paper is a formula for the structure sheaf
$\O_{\Omega_r}$ of a quiver variety in the Grothendieck ring $K^\circ
X$ of algebraic vector bundles on $X$.  This corresponds to computing
the alternating sum of a locally free resolution, so the $K$-theory
formula contains the cohomology formula as its leading term.

Our formula has the form
\begin{equation}
\label{eqn:K_quiver}
[\O_{\Omega_r}] = \sum_{|\mu| \geq d(r)} c_\mu(r)\, G_{\mu_1}(E_1 -
E_0) \cdots G_{\mu_n}(E_n - E_{n-1}) \in
K^\circ X
\end{equation}
where the sum is this time over a finite collection of sequences of
partitions for which the weights add up to at least the expected
codimension.  The elements \linebreak $G_{\mu_i}(E_i - E_{i-1}) \in
K^\circ X$ are called stable Grothendieck polynomials; these will be
defined in \refsec{sec:groth}.  The coefficients $c_\mu(r)$ in this
formula are given by a generalization of the algorithm for the
coefficients of (\ref{eqn:cohom_quiver}).  In particular, the
coefficients are the same when $|\mu| = d(r)$.

We conjecture that the signs of the new coefficients alternate with
the weight of $\mu$, {i.e.\ }$(-1)^{|\mu| - d(r)} c_\mu(r) \geq 0$.  It
appears to be a rather general phenomenon that coefficients which show
up in $K$-theoretic formulas tend to have alternating signs, although
this is very poorly understood.  For example, a formula of Fomin and
Kirillov shows that the signs of the coefficients in Grothendieck
polynomials alternate with degree \cite{fomin.kirillov:yang-baxter}.
Similarly we have proved in \cite{buch:littlewood-richardson} that the
structure constants of the Grothendieck ring of a Grassmann variety
with respect to its basis of Schubert structure sheaves have signs
which alternate with codimension.  In fact, this is a special case of
our conjecture, since said structure constants are special cases of
the coefficients $c_\mu(r)$ of (\ref{eqn:K_quiver}).  It is worth
pointing out that in all cases where alternation of signs in
$K$-theory has been proved, this has been achieved by giving explicit
formulas for the coefficients in question.  This is in contrast to
cohomology, where positivity results can often be obtained by
realizing coefficients as the number of points in an intersection of
varieties in general position.

Our conjecture is true when $\Omega_r$ is a variety of complexes,
{i.e.\ }$r_{ij} = 0$ whenever $j - i \geq 2$.  In fact, the algorithm
for the coefficients $c_\mu(r)$ is particularly simple in this case,
and it shows that (\ref{eqn:K_quiver}) is multiplicity free in the
sense that every coefficient $c_\mu(r)$ is either $1$, $-1$, or zero.
We have furthermore verified the conjecture computationally for all 
sequences with at most 4 bundles of ranks up to 7.

The proof of the cohomology formula in \cite{buch.fulton:chern} is
based on the simple idea of realizing the quiver variety $\Omega_r$ as
a birational image of a simpler quiver variety $\Omega_{\bar r}$ which
lives on a product of Grassmann bundles over $X$.  The class of
$\Omega_r$ can then be calculated inductively as the pushforward of
the class of $\Omega_{\bar r}$, which is done using a Gysin formula of
Pragacz \cite{pragacz:enumerative}.  However, before this Gysin
formula can be applied, one must first rearrange the inductive formula
for $\Omega_{\bar r}$ by replacing the Schur polynomials
$s_{\mu_i}(E_i - E_{i-1})$ in this formula with linear combinations of
products $s_\sigma(E_i - F) \cdot s_\tau(F - E_{i-1})$ for other
bundles $F$, which can be done by invoking the coproduct in the ring
of symmetric functions.  Thus the cohomology formula is a consequence
of the large cohomological toolbox surrounding the ring of symmetric
functions, once the right overall geometric construction has been
made.  Of particular importance are the coproduct on Schur functions
and Pragacz's Gysin formula, as well as the Thom-Porteous formula for
starting the induction.

While the same method turns out to work for the $K$-theory formula, it
was far from obvious that this would be possible when we started our
project.  First of all, while double stable Grothendieck polynomials
had been defined by Fomin and Kirillov
\cite{fomin.kirillov:yang-baxter, fomin.kirillov:grothendieck} and
studied combinatorially \cite{fomin.greene:noncommutative}, they had
never been applied to geometry.  Furthermore, the properties of Schur
functions that are needed for the cohomology formula had no known
analogues.  Our work on generalizing the formula has therefore
consisted mainly of finding and proving $K$-theoretic analogues of
known cohomological tools.


The first step in this direction was carried out in
\cite{buch:littlewood-richardson} where we proved that the linear span
of all stable Grothendieck polynomials form a bialgebra $\Gamma$ which
is a $K$-theory parallel of the ring of symmetric functions.  In the
same way as the ring of symmetric functions describes cohomology of
Grassmannians, $\Gamma$ describes their $K$-theory.  In this paper we
prove a $K$-theory version of the Thom-Porteous formula, a Gysin
formula for calculating $K$-theoretic pushforwards from a Grassmann
bundle which generalizes Pragacz's cohomological formula, and we
develop the few extra bits of combinatorics which make it all fit
together.

One additional ingredient in the proofs of (\ref{eqn:cohom_quiver})
and (\ref{eqn:K_quiver}) is a result of Lakshmibai and Magyar showing
that a quiver variety of the expected codimension is Cohen-Macaulay
\cite{lakshmibai.magyar:degeneracy}.  For the $K$-theory formula we
shall furthermore need their result about rational singularities of
quiver varieties to deduce that the structure sheaf of the inductive
quiver variety $\Omega_{\bar r}$ mentioned above pushes forward to the
structure sheaf of $\Omega_r$.


In \refsec{sec:groth} we fix the notation regarding Grothendieck
polynomials and stable Grothendieck polynomials, and we explain their
relations to geometry.  In \refsec{sec:intseq} we define stable
Grothendieck polynomials for arbitrary sequences of integers which
extend the definition of stable Grothendieck polynomials for
partitions.  This is needed for describing the algorithm for the
coefficients in (\ref{eqn:K_quiver}).  This algorithm is then
presented in \refsec{sec:formula}, where we also explain the meaning
of our formula when $X$ is singular or $\Omega_r$ does not have its
expected codimension.  In addition we interpret the formula in the
case of varieties of complexes.  In \refsec{sec:appgroth} we show that
Grothendieck polynomials and stable Grothendieck polynomials are
special cases of the quiver formula.  Combined with some recent
results of Lascoux \cite{lascoux:transition}, this supplies additional
evidence for our conjecture about the signs of the coefficients
$c_\mu(r)$.  The last two sections are devoted to proving our
generalization of Pragacz's Gysin formula.  \refsec{sec:jacobi_trudi}
proves a generalization of the Jacobi-Trudi formula for Schur
functions, which in \refsec{sec:gysin} is used to establish the Gysin
formula itself.  We finish the paper by noticing that the pushforward
map from a Grassmann bundle is multiplicative when applied to products
of Grothendieck polynomials for short partitions.

We are indebted to S.~Fomin for sharing some important insights in the
start of this project, which significantly improved our understanding
of stable Grothendieck polynomials.  In addition we thank W.~Fulton
for numerous helpful comments, suggestions, and encouragements during
the project.


%% file: groth.tex
\section{Grothendieck polynomials}
\label{sec:groth}

In this section we fix the notation concerning Grothendieck
polynomials and stable Grothendieck polynomials, and explain their
relations to geometry.  We furthermore summarize the necessary results
from \cite{buch:littlewood-richardson}.

Given a permutation $w \in S_n$, Lascoux and Sch{\"u}tzenberger define
the {\em double Grothendieck polynomial\/} $\Groth_w = \Groth_w(x;y)$
for $w$ as follows \cite{lascoux.schutzenberger:structure}.  For the
longest permutation $w_0 = n \, (n-1) \cdots 2 \, 1$ we set
\[ \Groth_{w_0} = \prod_{i+j \leq n} (x_i + y_j - x_i y_j) \,. \]
If $w$ is not the longest permutation, we can find a simple reflection
$s_i = (i,i+1) \in S_n$ such that $\ell(w s_i) = \ell(w) + 1$.  Here
$\ell(w)$ denotes the length of $w$, which is the smallest number
$\ell$ for which $w$ can be written as a product of $\ell$ simple
reflections.  We then define
\[ \Groth_w = \pi_i(\Groth_{w s_i}) \]
where $\pi_i$ is the isobaric divided difference operator given by
\[ \pi_i(f) = \frac{(1-x_{i+1}) f(x_1,x_2,\dots) - (1-x_i)
  f(\dots,x_{i+1},x_i,\dots)}{x_i - x_{i+1}} \,. 
\]
This definition is independent of our choice of simple reflection
$s_i$ since the operators $\pi_i$ satisfy the Coxeter relations.


Notice that the longest element in $S_{n+1}$ is $w_0^{(n+1)} = w_0
\cdot s_n \cdot s_{n-1} \cdots s_1$.  Since $\pi_n \cdot \pi_{n-1}
\cdots \pi_1$ applied to the Grothendieck polynomial for $w_0^{(n+1)}$
is equal to $\Groth_{w_0}$, it follows that $\Groth_w$ does not depend
on which symmetric group $w$ is considered an element of.

Now let $F_1 \subset F_2 \subset \dots \subset F_n
\xrightarrow{\varphi} H_n \twoheadrightarrow \dots
\twoheadrightarrow H_2 \twoheadrightarrow H_1$ be a full flag of
vector bundles on $X$ followed by a map $\varphi$ to a dual full flag.
For $w \in S_{n+1}$ we define the degeneracy locus
\[ \Omega_w = \Omega_w(F_\bull \to H_\bull) = 
   \{ x \in X \mid \rank(F_q(x) \to H_p(x)) \leq r_w(p,q) 
   ~\forall p,q \}
\]
where $r_w(p,q) = \# \{ i \leq p \mid w(i) \leq q \}$.  The expected
codimension for this locus is the length of $w$.

Suppose $F_1 \subset F_2 \subset \dots \subset F_n \subset V$ is a
full flag of subbundles in a vector bundle $V$ of rank $n+1$.  Let
$\pi : \Fl^*(V) \to X$ be the bundle of dual flags of $V$, with
tautological flag $\pi^* V \twoheadrightarrow \Tilde H_n
\twoheadrightarrow \dots \twoheadrightarrow \Tilde H_1$.  In this case
the Schubert variety $\Tilde \Omega_w = \Omega_w(\pi^* F_\bull \to
\Tilde H_\bull)$ has codimension $\ell(w)$ in $\Fl^*(V)$.  Fulton and
Lascoux \cite{fulton.lascoux:pieri} have proved that its structure
sheaf is given by the double Grothendieck polynomial for $w$:
\begin{equation}
\label{eqn:fl}
  [\O_{\Tilde \Omega_w}] = \Groth_w
  (1-\Tilde L_1^{-1}, \dots, 1-\Tilde L_n^{-1}; 1-M_1,\dots,1-M_n)
\end{equation}
in $K^\circ \Fl^*(V)$, where $\Tilde L_i = \ker(\Tilde H_i \to \Tilde
H_{i-1})$ and $M_i = F_i/F_{i-1}$.

Using the fact that a Grothendieck polynomial $\Groth_w(x;y)$ does not
depend on which symmetric group $w$ belongs to, this formula readily
generalizes as follows:

\begin{thm} 
\label{thm:doubleflag}
  If the locus $\Omega_w = \Omega_w(F_\bull \to H_\bull)$ in $X$ has
  its expected codimension $\ell(w)$, then 
\[ [\O_{\Omega_w}] = 
  \Groth_w(1 - L_1^{-1}, \dots, 1 - L_n^{-1}; 1-M_1, \dots, 1-M_n)
\]
where $L_i = \ker(H_i \to H_{i-1})$ and $M_i = F_i/F_{i-1}$.
\end{thm}
\begin{proof}
  Set $V = F_n \oplus H_n$ and let $\pi : \Fl^*(V) \to X$ be the dual
  flag bundle of $V$ with tautological flag $\pi^* V
  \twoheadrightarrow \Tilde H_{2n-1} \twoheadrightarrow \dots
  \twoheadrightarrow \Tilde H_1$.  Define $\psi : F_n \to V$ by
  $\psi(\sigma) = (\sigma, \varphi(\sigma))$, and set $F_i = \psi(F_n)
  + \ker(H_n \to H_{2n-i})$ for $n < i < 2n$.  Then $F_1 \subset \dots
  \subset F_{2n-1} \subset V$ is a full flag of subbundles in $V$, and
  (\ref{eqn:fl}) applies to give a formula for the structure sheaf of
  the Schubert variety $\Tilde \Omega_w = \Omega_w(\pi^* F_\bull \to
  \Tilde H_\bull)$.
  
  Set $H_i = (F_n/F_{2n-i}) \oplus H_n$ for $n < i < 2n$.  Then there
  is a unique section $s : X \to \Fl^*(V)$ such that the dual flag $V
  \twoheadrightarrow H_{2n-1} \twoheadrightarrow \dots
  \twoheadrightarrow H_1$ is the pullback of the tautological flag
  $\pi^* V \twoheadrightarrow \Tilde H_\bull$ on $\Fl^*(V)$, and
  furthermore we have $\Omega_w = s^{-1}(\Tilde \Omega_w)$ as
  subschemes of $X$.  Since the loci $\Omega_w$ and $\Tilde \Omega_w$
  have the same codimensions and are Cohen-Macaulay, this implies that
\[ \begin{split} 
  [\O_{\Omega_w}] &= s^* [\O_{\Tilde \Omega_w}] = 
  s^*\, \Groth_w(1-\Tilde L_1^{-1}, \dots, 1-\Tilde L_{2n-1}^{-1}; 
               1-M_1, \dots, 1-M_{2n-1}) \\
  &= s^*\, \Groth_w(1-\Tilde L_1^{-1}, \dots, 1-\Tilde L_n^{-1}; 
                  1-M_1, \dots, 1-M_n) \\
  &= \Groth_w(1-L_1^{-1}, \dots, 1-L_n^{-1}; 1-M_1, \dots, 1-M_n)
\end{split} \]
which completes the proof.
\end{proof}

We now turn to stable Grothendieck polynomials.  Given a permutation
$w \in S_n$ and a non-negative integer $m$, we let $1^m \times w \in
S_{m+n}$ denote the shifted permutation which is the identity on $\{1,
2, \dots, m\}$ and which maps $j$ to $w(j-m) + m$ for $j > m$.  Fomin
and Kirillov have shown that when $m$ grows to infinity, the
coefficient of each fixed monomial in $\Groth_{1^m \times w}$
eventually becomes stable \cite{fomin.kirillov:yang-baxter}.  The {\em
  double stable Grothendieck polynomial\/} $G_w \in \Z\llbracket
x_i,y_i \rrbracket_{i \geq 1}$ is defined as the resulting power
series:
\[ G_w = G_w(x;y) = \lim_{m \to \infty} \Groth_{1^m \times w} \,. \]
Fomin and Kirillov also proved that this power series is symmetric in
the variables $\{x_i\}$ and $\{y_i\}$ separately, and that
\[ G_w(1-e^{-x}; 1-e^y) = 
   G_w(1-e^{-x_1}, 1-e^{-x_2}, \dots; 1-e^{y_1}, 1-e^{y_2}, \dots) 
\]
is super symmetric, i.e.\ if one sets $x_1 = y_1$ in this expression
then the result is independent of $x_1$ and $y_1$.  Alternatively,
these facts can be deduced from \refthm{thm:doubleflag}.

We shall be mostly concerned with stable Grothendieck polynomials for
Grassmannian permutations.  If $\lambda$ is a partition and $p \geq
\ell(\lambda)$, {i.e.\ }$\lambda_{p+1} = 0$, the {\em Grassmannian
  permutation\/} for $\lambda$ with descent in position $p$ is the
unique permutation $w_\lambda$ such that $w_\lambda(i) = i +
\lambda_{p+1-i}$ for $1 \leq i \leq p$ and $w_\lambda(i) <
w_\lambda(i+1)$ for $i \neq p$.  We define $G_\lambda =
G_{w_\lambda}$.  Notice that if $q > p$, then the Grassmannian
permutation for $\lambda$ with descent at position $q$ is equal to
$1^{q-p} \times w_\lambda$.  Therefore $G_\lambda$ is independent of
the choice of $p$.

Let $\Gamma \subset \Z \llbracket x_i, y_i \rrbracket$ be the linear
span of all stable Grothendieck polynomials.  It is shown in
\cite{buch:littlewood-richardson} that this group is a bialgebra and
that the elements $G_\lambda$ form a basis.  We will proceed to
describe the structure constants of $\Gamma$.

If $a$ and $b$ are two non-empty subsets of the positive integers
$\N$, we will write $a < b$ if $\max(a) < \min(b)$, and $a \leq b$ if
$\max(a) \leq \min(b)$.  We define a {\em set-valued tableau\/} to be
a labeling of the boxes in a Young diagram or skew diagram with finite
non-empty subsets of $\N$, such that the rows are weakly increasing
from left to right and the columns strictly increasing from top to
bottom.  The shape $\sh(T)$ of a tableau $T$ is the partition or skew
diagram it is a labeling of.  For example,
\[ \tableau{ & & {1} & {2\,3} \\ & {1\,2} & {234} \\
   {2} & {3\,5} & {7}}
\]
is a set-valued tableau whose shape is the skew diagram between the
partitions $(4,3,3)$ and $(2,1)$.  The {\em word\/} of a set-valued
tableau is the sequence of integers in its boxes when these are read
left to right and then bottom to top, and the integers in a single box
are arranged in increasing order.  The word of the above tableau is
$(2,3,5,7,1,2,2,3,4,1,2,3)$.

We say that a sequence of positive integers $w = (i_1, i_2, \dots,
i_\ell)$ is has {\em content\/} $(c_1, c_2, \dots, c_r)$ if $w$
consists of $c_1$ 1's, $c_2$ 2's, and so on up to $c_r$ $r$'s.  If the
content of each subsequence $(i_k, \dots, i_\ell)$ of $w$ is a
partition, then $w$ is called a {\em reverse lattice word}.

If $\lambda$ and $\mu$ are partitions, we let $\lambda * \mu$ denote
the skew diagram obtained by attaching the Young diagrams for
$\lambda$ and $\mu$ corner to corner as shown:
\[ \lambda * \mu ~=~ 
   \raisebox{-32pt}{\includegraphics[scale=0.7]{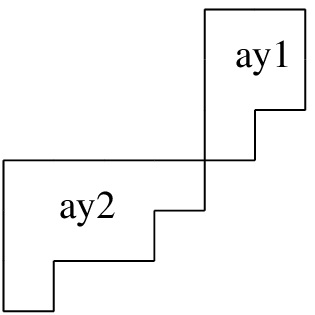}} 
\]
The main result of \cite{buch:littlewood-richardson} now says that 
\begin{equation}
\label{eqn:lrmult}
  G_\lambda \cdot G_\mu = \sum_\nu c^\nu_{\lambda \mu} \, G_\nu
\end{equation}
where $c^\nu_{\lambda \mu}$ is equal to $(-1)^{|\nu|-|\lambda|-|\mu|}$
times the number of set-valued tableaux $T$ of shape $\lambda * \mu$
such that the word of $T$ is a reverse lattice word with content
$\nu$.  

Now if $\lambda$, $\mu$, and $\nu$ are partitions, we set
$d^\nu_{\lambda \mu} = c^\rho_{\nu R}$, where $R = (p)^q$ is any
rectangular partitions containing $\lambda$ and $\mu$, and $\rho =
(p+\lambda_1, \dots, p+\lambda_q,\mu_1, \mu_2, \dots)$ is the
partition obtained by attaching $\lambda$ and $\mu$ to the sides of
$R$.  In \cite{buch:littlewood-richardson} it is proved that these
coefficients do not depend on the choice of the rectangle $R$.
Furthermore, whenever $x$, $y$, $z$, and $w$ are different sets of
variables we have
\begin{equation}
\label{eqn:coprod}
  G_\nu(x,y; z,w) = \sum_{\lambda,\mu} d^\nu_{\lambda \mu} \,
  G_\lambda(x;z) \cdot G_\mu(y;w) \,.
\end{equation}

\begin{thm}
\label{thm:bialg}
  The group $\Gamma = \bigoplus_\lambda G_\lambda \subset \Z
  \llbracket x_i, y_i \rrbracket$ is a commutative and cocommutative
  bialgebra with unit and counit.  Multiplication is given by
  (\ref{eqn:lrmult}) and the coproduct $\Delta : \Gamma \to \Gamma
  \otimes \Gamma$ is defined by $\Delta G_\nu = \sum_{\lambda,\mu}
  d^\nu_{\lambda \mu} G_\lambda \otimes G_\mu$.
\end{thm}

It is also possible to give a formula for stable Grothendieck
polynomials based on set-valued tableaux.  Given a tableau $T$, let
$x^T$ be the monomial in which the exponent of $x_i$ is the number of
boxes in $T$ which contain the integer $i$.  If $T$ is the tableau
displayed above we get $x^T = x_1^2\, x_2^4\, x_3^3\, x_4\, x_5 \,
x_7$.  We let $|T|$ denote the total degree of this monomial, i.e.\ 
the sum of the cardinalities of the sets in the boxes of $T$.  In 
\cite{buch:littlewood-richardson} it is proved that the single stable Grothendieck polynomial $G_\lambda(x) = G_\lambda(x;0)$ is given by
\begin{equation}
\label{eqn:single}
 G_\lambda(x) = \sum_{\sh(T) = \lambda} (-1)^{|T|-|\lambda|} \, x^T 
\end{equation}
where the sum is over all set-valued tableaux $T$ of shape $\lambda$.
The double stable Grothendieck polynomial for $\lambda$ is then given
by
\[ G_\lambda(x;y) = \sum_{\sigma,\tau} d^\lambda_{\sigma \tau} \,
   G_\sigma(x) \cdot G_{\tau'}(y)
\]
where $\tau'$ denotes the conjugate partition of $\tau$.

Let $F = L_1 \oplus \dots \oplus L_f$ and $E = M_1 \oplus \dots \oplus
M_e$ be vector bundles on a variety $X$ which are both direct sums of
line bundles.  We then define
\[ G_\lambda(F - E) = G_\lambda(1-L_1^{-1}, \dots, 1-L_f^{-1};
   1-M_1,\dots,1-M_e) ~\in~ K^\circ(X) \,.
\]
Since $G_\lambda$ is symmetric, this is a polynomial in the exterior
powers of $F^\vee$ and $E$.  Therefore the definition makes sense also
when $E$ and $F$ are not direct sums of line bundles.  For example we
have $G_1(F - E) = 1 - \frac{\bigwedge^e E}{\bigwedge^f F}$.  The fact
that $G_\lambda(1-e^{-x}; 1-e^y)$ is super-symmetric implies that
$G_\lambda(F \oplus H - E \oplus H) = G_\lambda(F - E)$ for any bundle
$H$.  In particular we can regard $G_\lambda$ as a well defined
function $G_\lambda : K^\circ X \to K^\circ X$.  Equation
(\ref{eqn:coprod}) then says that for any elements $\alpha, \beta \in
K^\circ X$ we have
\begin{equation} \label{eqn:split}
G_\lambda(\alpha + \beta) = \sum_{\sigma,\tau}
d^\lambda_{\sigma \tau} \, G_\sigma(\alpha) \cdot G_\tau(\beta) \,.
\end{equation}
Another useful fact, due to Fomin, is that $G_\lambda(F - E) =
G_{\lambda'}(E^\vee - F^\vee)$.

This notation makes it possible to give a Thom-Porteous formula for
$K$-theory which is analogous to its cohomological equivalent.  Let $E
\to F$ be a morphism between vector bundles of ranks $e$ and $f$.
Given an integer $r \leq \min(e,f)$ we have the degeneracy locus
$\Omega_r(E \to F) = \{ x \in X \mid \rank(E(x) \to F(x)) \leq r \}$.

\begin{thm}
\label{thm:porteous}
If the codimension of $\Omega_r(E \to F)$ in $X$ is $(e-r)(f-r)$ then
the class of its structure sheaf is given by
\[ [\O_{\Omega_r(E \to F)}] = G_\lambda(F - E) \]
where $\lambda = (e-r)^{f-r}$ is a rectangular partition with $f-r$
rows and $e-r$ columns.
\end{thm}
\begin{proof}
By the splitting principle we may assume that $E$ and $F$ come
equipped with full flags $E_1 \subset \dots \subset E_e = E$ and $F =
F_f \twoheadrightarrow \dots \twoheadrightarrow F_1$.  Let $w_\lambda$
be the Grassmannian permutation for $\lambda$ with descent at position
$f$.  Then $w$ is a permutation in $S_n$ where $n = e+f-r$.  Set $E_i
= E \oplus \O_X^{\oplus i-e}$ for $e < i < n$ and $F_j = F \oplus
\O_X^{\oplus j-f}$ for $f < j < n$, and let $\varphi : E_{n-1} \to
F_{n-1}$ be the map $E_{n-1} \to E \to F \to F_{n-1}$, i.e.\ the map
$E \to F$ is extended by zeros on the trivial parts of $E_{n-1}$ and
$F_{n-1}$.  It is now easy to check that $\Omega_r(E \to F) =
\Omega_{w_\lambda}(E_\bull \to F_\bull)$ as subschemes of $X$, so by
\refthm{thm:doubleflag} we get
\[ [\O_{\Omega_r(E \to F)}] = \Groth_{w_\lambda}(1-L_1^{-1}, \dots,
   1-L_f^{-1}, 0, \dots, 0; 1-M_1, \dots, 1-M_e, 0, \dots, 0)
\]
where $L_j = \ker(F_j \to F_{j-1})$ and $M_i = E_i/E_{i-1}$.  Notice
that $\Omega_{w_\lambda}(E_\bull \to F_\bull) = \Omega_{1 \times
  w_\lambda}(E_\bull \subset E_{n-1} \oplus \O_X \xrightarrow{\varphi
  \oplus 1} F_{n-1} \oplus \O_X \twoheadrightarrow F_\bull)$.  This
means that the formula does not change when we shift the permutation
$w_\lambda$, so in fact we have
\[ [\O_{\Omega_r(E \to F)}] = G_{w_\lambda}(1-L_1^{-1}, \dots,
   1-L_f^{-1}; 1-M_1, \dots, 1-M_e) = G_\lambda(F - E) \,.
\]
This finishes the proof.
\end{proof}


%% file: intseq.tex
\section{Sequences of integers}
\label{sec:intseq}

\newcommand{\GG}{{\mathcal G}}

In order to define the coefficients in our formula for quiver
varieties, we need to define stable Grothendieck polynomials for
arbitrary sequences of integers.  Our definition of these is inspired
by the following determinant formula of Lenart.

For integers $k \in \Z$ and $i \geq 0$, let $h_k(x_1,\dots,x_n/1^i)$
denote the coefficient of $t^k$ in the formal power series expansion
of
\[ \frac{(1-t)^i}{\prod_{j=1}^n (1-x_j t)} \,. \]
In particular we have $h_0(x_1,\dots,x_n/1^i) = 1$ and
$h_k(x_1,\dots,x_n/1^i) = 0$ for $k < 0$.

Let $I = (I_1, I_2, \dots, I_\ell)$ be a finite sequence of integers
of length $\ell$.  For convenience we shall regard $I_j$ as being zero
if $j > \ell$.  For $n > \ell$ we now define $\GG_I(x_1,\dots,x_n)$ to
be the determinant of the $n \times n$ matrix whose $(i,j)$'th entry
is equal to $h_{I_i+j-1}(x_1,\dots,x_n/1^{i-1})$:
\[ \GG_I(x_1,\dots,x_n) = \det \big( 
   h_{I_i+j-1}(x_1,\dots,x_n/1^{i-1}) \big)_{1 \leq i,j \leq n}
\]
Notice that the size of this determinant depends on the number of
variables.  With this notation we have:

\begin{thm}[Lenart \cite{lenart:combinatorial}]
  If $I$ is a partition then
\[ \GG_I(x_1,\dots,x_n) = \Groth_{w_I}(x) = G_I(x_1,\dots,x_n) \,. \]
\end{thm}

\begin{lemma}
\label{lemma:permute}
Let $I$ and $J$ be sequences of integers and suppose $p < q$ are
integers.  Then
\[ \GG_{I,p,q,J}(x_1,\dots,x_n) = 
   \sum_{k=p+1}^q \GG_{I,q,k,J}(x_1,\dots,x_n) - 
   \sum_{k=p+1}^{q-1} \GG_{I,q-1,k,J}(x_1,\dots,x_n) \,.
\]
\end{lemma}
\begin{proof}
To cut down on the notation, we shall prove this in the case where $I$
and $J$ are empty and $n=2$.  The proof of the general case is exactly
the same.  For convenience we will also write $h_k(x/1^i)$ for
$h_k(x_1,x_2/1^i)$.

Using the rule $h_k(x) = h_{k+1}(x) - h_{k+1}(x/1)$
repeatedly we get
\[
\begin{vmatrix}
h_p(x) & h_{p+1}(x) \\
h_q(x) & h_{q+1}(x)
\end{vmatrix} =
\sum_{k=p+1}^q
\begin{vmatrix}
-h_k(x/1) & -h_{k+1}(x/1) \\
h_q(x) & h_{q+1}(x)
\end{vmatrix} 
= \sum_{k=p+1}^q \GG_{q,k}(x_1,x_2) \,.
\]
The lemma follows from this since 
\[ 
\begin{vmatrix}
h_p(x) & h_{p+1}(x) \\
h_q(x/1) & h_{q+1}(x/1)
\end{vmatrix} 
=
\begin{vmatrix}
h_p(x) & h_{p+1}(x) \\
h_q(x) & h_{q+1}(x)
\end{vmatrix} 
-
\begin{vmatrix}
h_p(x) & h_{p+1}(x) \\
h_{q-1}(x) & h_q(x)
\end{vmatrix} 
\,.
\]
\end{proof}

\begin{cor}
  Let $I = (I_1, I_2, \dots, I_\ell)$ be a sequence of integers and
  let $n$ be an integer such that $n \geq \ell$ and $n \geq i - I_i$
  for all $1 \leq i \leq \ell$.  Then $\GG_I(x_1,\dots,x_n)$ is a
  finite linear combination of determinants
  $\GG_\lambda(x_1,\dots,x_n)$ for partitions $\lambda$:
\[ \GG_I(x_1,\dots,x_n) = \sum_\lambda \delta_{I,\lambda} \,
   \GG_\lambda(x_1,\dots,x_n) \,. 
\]
Furthermore the coefficients $\delta_{I,\lambda}$ do not depend on
$n$.
\end{cor}
\begin{proof}
  Define a ``potential'' function $\rho(I) = \sum_{j=1}^n (n - j)
  (\Bar I_j - I_j)$ where $\Bar I_j = \max \{ I_k \mid j \leq k \leq n
  \}$.  Then $\rho(I) \geq 0$ for all sequences $I$.
  
  We proceed by induction on $\rho(I)$.  If $\rho(I) = 0$ then $I$
  must be weakly decreasing.  In fact it must be a partition because
  the assumption $n \geq n - I_n$ implies that $I_n \geq 0$.
  Therefore $\GG_I(x_1,\dots,x_n)$ already has the desired form.
  
  If $\rho(I) > 0$ then for some $1 \leq j < n$ we must have $I_j <
  I_{j+1}$.  We can now apply \reflemma{lemma:permute} with $p = I_j$
  and $q = I_{j+1}$ to write $\GG_I(x_1,\dots,x_n)$ as a linear
  combination of other determinants $\GG_J(x_1,\dots,x_n)$, and it is
  easy to check that these satisfy $\rho(J) < \rho(I)$ and $n \geq i -
  J_i$ for all $1 \leq i \leq \ell$.  Each of these new determinants
  is therefore a linear combination of the polynomials
  $\GG_\lambda(x_1,\dots,x_n)$ by induction, which proves the claim
  for the sequence $I$.

  The fact that the coefficients $\delta_{I,\lambda}$ are independent
  of $n$ follows because the formula of \reflemma{lemma:permute} is
  independent of $n$.
\end{proof}

Now define $G_I = G_I(x;y) = \sum_\lambda \delta_{I,\lambda} \,
G_\lambda(x;y) \in \Gamma$.  This is well defined by the corollary,
and since $G_I(x_1,\dots,x_n) = \GG_I(x_1,\dots,x_n)$ when $n$ is
sufficiently large, we have $G_I(x) = \lim_{n \to \infty}
\GG_I(x_1,\dots,x_n)$.  Furthermore, \reflemma{lemma:permute} implies
that
\begin{equation}
\label{eqn:intseq}
  G_{I,p,q,J} = \sum_{k=p+1}^q G_{I,q,k,J} - 
  \sum_{k=p+1}^{q-1} G_{I,q-1,k,J} 
\end{equation}
whenever $p < q$.  This gives a practical way to compute the
polynomials $G_I$.

\begin{example}
$G_{1\,1\,3} = G_{1\,3\,2} + G_{1\,3\,3} - G_{1\,2\,2} = G_{3\,2\,2} + G_{3\,3\,2} + G_{3\,3\,3} - 2\,G_{2\,2\,2}$.
\end{example}

Contrary to the case of Schur functions, a stable Grothendieck
polynomial $G_I$ for a sequence of integers is never equal to zero.
In fact, (\ref{eqn:intseq}) readily implies that $\sum_\lambda
\delta_{I,\lambda} = 1$ for any sequence of integers $I$.  It is also
easy to prove that if $J$ is a sequence of non-positive integers, then
$G_{I,J} = G_I$.  In addition, if $G_\lambda$ occurs in the expansion
of $G_I$, then $\lambda$ must be contained in the partition $\Bar I =
(\Bar I_1, \Bar I_2, \dots)$, and furthermore we have $\delta_{I,\Bar
  I} = 1$.  A lower bound on $\lambda$ may also be obtained.  Let
$\rho = (0,1,2,\dots)$ and let $J$ denote the sequence $I - \rho =
(I_1, I_2 - 1, I_3 - 2, \dots)$ arranged in decreasing order.  Then
any partition $\lambda$ for which $G_\lambda$ occurs in $G_I$ must
contain the partition $\Tilde I = \overline{J+\rho}$.  This lower
bound is not sharp.  If we take $I = (0,2,0,3)$ then $\Tilde I =
(2,2,2,1)$ and $\delta_{I,\Tilde I} = 0$.  We will not need these
remarks in the following.


%% file: formula.tex
\section{A formula for quiver varieties}
\label{sec:formula}

We are now ready to describe our formula for the structure sheaf of a
quiver variety.  Let $X$ be any Noetherian scheme equipped with a
sequence $E_\bull$ of vector bundles, and let $\Omega_r =
\Omega_r(E_\bull)$ be the associated quiver variety.  We define a
localized class $\OOmega_r$ in the Grothendieck group $K_\circ
\Omega_r$ of coherent sheaves on $\Omega_r$ as follows.  On the bundle
$H = \Hom(E_0, E_1) \times_X \dots \times_X \Hom(E_{n-1}, E_n)
\xrightarrow{~\pi~} X$ we have a sequence of tautological maps $\pi^*
E_0 \to \pi^* E_1 \to \dots \to \pi^*E_n$.  We let $\Tilde \Omega_r
\subset H$ denote the quiver variety defined by this sequence.  Now
the bundle maps on $X$ define a section $s : X \to H$, and $\Omega_r =
s^{-1}(\Tilde \Omega_r)$.  The localized class $\OOmega_r$ is defined
by
\begin{equation}
\label{eqn:class}
  \OOmega_r = s^!( [\O_{\Tilde \Omega_r}] ) = 
  \sum_{j \geq 0} (-1)^j \Tor_j^{H}(\O_X, \O_{\Tilde \Omega_r}) 
  \in K_\circ \Omega_r \,.
\end{equation}
Notice that since $s$ is a regular embedding, it follows that locally
on $H$ the structure sheaf of $X$ has a finite free resolution, so the
sum in (\ref{eqn:class}) is finite.  The definition of $\OOmega_r$
implies that these classes are compatible with perfect pullback and
proper pushforward \cite{fulton.macpherson:categorical}.

The codimension of $\Tilde \Omega_r$ in $H$ is always equal to $d(r)$
\cite{buch.fulton:chern}.  Furthermore, Lakshmibai and Magyar have
shown that this locus is Cohen-Macaulay and has rational singularities
if $X$ has these properties \cite{lakshmibai.magyar:degeneracy}.  If
$X$ is Cohen-Macaulay and $\Omega_r$ has its expected codimension
$d(r)$ in $X$, this implies that $\Omega_r$ is Cohen-Macaulay as well.
In addition, a local regular sequence generating the ideal of $X$ in
$H$ pulls back to a local regular sequence defining the ideal of
$\Omega_r$ in $\Tilde \Omega_r$ \cite[Lemma
A.7.1]{fulton:intersection}.  It follows from this that
$\Tor^H_j(\O_X, \O_{\Tilde \Omega_r}) = 0$ for all $j > 0$, so
$\OOmega_r = [\O_X \otimes_{\O_H} \O_{\Tilde \Omega_r}] =
[\O_{\Omega_r}]$.  More generally, this is true without the
Cohen-Macaulay condition if $\depth(\Omega_r, X) = d(r)$
\cite[Ex.~14.3.1]{fulton:intersection}.  We can now state our main
result.

\begin{thm}
\label{thm:structure}
The image of $\OOmega_r$ in $K_\circ X$ is given by
\[
  \OOmega_r = \sum_{|\mu| \geq d(r)} c_\mu(r) \, 
  G_{\mu_1}(E_1-E_0) \cdots G_{\mu_n}(E_n - E_{n-1}) \,.
\]
The sum is over a finite number of sequences $\mu$ of partitions
$\mu_i$ such that the sum of the weights of these partitions is at
least equal to $d(r)$.  The coefficients $c_\mu(r)$ are integers which
are given by an explicit combinatorial algorithm.
\end{thm}

The algorithm which computes the coefficients $c_\mu(r)$ is the same
as the one computing the coefficients in the cohomology formula
\cite{buch.fulton:chern}, except the bialgebra $\Gamma$ replaces the
ring of symmetric functions.  To describe the algorithm, we will
construct an element $P_r$ in the $n$th tensor power of $\Gamma$, such
that
\[ P_r = \sum_\mu c_\mu(r) \, 
   G_{\mu_1} \otimes \dots \otimes G_{\mu_n} \,. 
\]

It is convenient to arrange the rank conditions in a {\em rank diagram\/}:
\[ \begin{matrix}
E_0 & \to & E_1 & \to & E_2 & \to & \cdots & \to & E_n
\vspace{0.1cm} \\
r_{00} && r_{11} && r_{22} && \cdots && r_{nn} \\
& r_{01} && r_{12} && \cdots && r_{n-1,n} \\
&& r_{02} && \cdots && r_{n-2,n} \\
&&& \ddots \\
&&&& r_{0n}
\end{matrix} \]
In this diagram we replace each small triangle of numbers
\[ \begin{matrix}
r_{i,j-1} && r_{i+1,j} \\
& r_{ij}
\end{matrix} \]
with a rectangle $R_{ij}$ with $r_{i+1,j} - r_{ij}$ rows and $r_{i,j-1}
- r_{ij}$ columns.
\[ R_{ij} = \raisebox{-18pt}{\picC{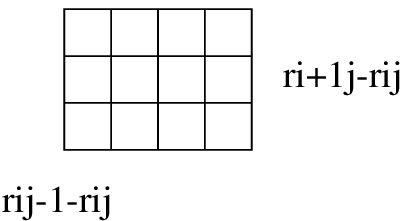}} \]
These rectangles are then arranged in a {\em rectangle diagram\/}:
\[ \begin{matrix}
R_{01} && R_{12} && \cdots && R_{n-1,n} \\
& R_{02} && \cdots && R_{n-2,n} \\
&& \ddots \\
&&& R_{0n}
\end{matrix} \]

The combinatorial data contained in the rank conditions $r = \{ r_{ij}
\}$ is very well represented by this diagram.  First of all, the rank
conditions can occur if and only if the rectangles always get shorter
when one travels south-east, while they get narrower when one travels
south-west.  Furthermore, the expected codimension $d(r)$ is equal to
the sum of the areas of the rectangles $R_{ij}$.  Finally, the element
$P_r$ depends only on the rectangle diagram.

We will define $P_r \in \Gamma^{\otimes n}$ by induction on $n$.  When
$n = 1$ (corresponding to a sequence of two vector bundles), the
rectangle diagram has only one rectangle $R = R_{01}$.  In this case
we set
\[ P_r = G_R \in \Gamma^{\otimes 1} \]
where $R$ is identified with the partition for which it is the Young
diagram.  This case recovers the Thom-Porteous formula
(\refthm{thm:porteous}).

If $n \geq 2$ we let $\Bar r$ denote the bottom $n$ rows of the rank
diagram.  Then $\Bar r$ is a valid set of rank conditions, so by
induction we can assume that 
\begin{equation}
\label{eqn:prbar}
P_{\Bar r} = \sum_\mu c_\mu(\Bar r) \, 
G_{\mu_1} \otimes \dots \otimes G_{\mu_{n-1}}
\end{equation}
is a well defined element of $\Gamma^{\otimes n-1}$.  Now $P_r$ is
obtained from $P_{\Bar r}$ by replacing each basis element $G_{\mu_1}
\otimes \dots \otimes G_{\mu_{n-1}}$ in (\ref{eqn:prbar}) with the sum
\[ \sum_{\stackrel{\sigma_1,\dots,\sigma_{n-1}}{\tau_1,\dots,\tau_{n-1}}}
   \left(\prod_{i=1}^{n-1} d^{\mu_i}_{\sigma_i \tau_i}\right)
   G_{\picC{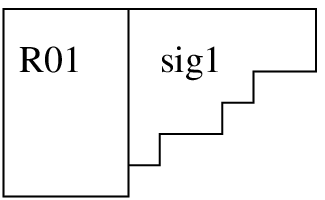}} \otimes \cdots \otimes
   G_{\picC{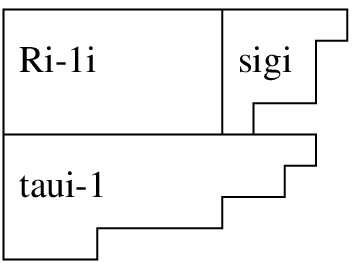}} \otimes \cdots \otimes
   G_{\picC{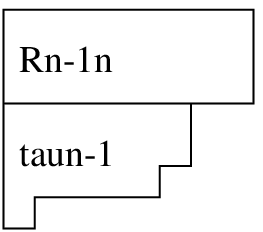}} \,.
\]
This sum is over all partitions $\sigma_1, \dots, \sigma_{n-1}$ and
$\tau_1, \dots, \tau_{n-1}$ such that $\sigma_i$ has fewer rows than
$R_{i-1,i}$ and the coproduct structure constant $d^{\mu_i}_{\sigma_i
  \tau_i}$ of $\Gamma$ is non-zero.  A diagram consisting of a
rectangle $R_{i-1,i}$ with (the Young diagram of) a partition
$\sigma_i$ attached to its right side, and $\tau_{i-1}$ attached
beneath should be interpreted as the sequence of integers giving the
number of boxes in each row of this diagram.

It can happen that the rectangle $R_{i-1,i}$ is empty, since the
number of rows or columns can be zero.  If the number of rows is zero,
then $\sigma_i$ is required to be empty, and the diagram is the Young
diagram of $\tau_{i-1}$.  If the number of columns is zero, then the
algorithm requires that the length of $\sigma_i$ is at most equal to
the number of rows $r_{ii} - r_{i-1,i}$ of $R_{i-1,i}$, and the
diagram consists of $\sigma_i$ in the top $r_{ii} - r_{i-1,i}$ rows
and $\tau_{i-1}$ below this, possibly with some zero-length rows in
between.

The proof of \refthm{thm:structure} is word for word identical to the
proof given in \cite{buch.fulton:chern}, except that the lemmas 2, 3,
and 4 of \cite{buch.fulton:chern} are replaced with \refthm{thm:gysin}
from this paper, Corollary 6.5 of \cite{buch:littlewood-richardson},
and equation (\ref{eqn:split}), respectively.  The only point which
requires a comment is that the modified proof will need that for
certain proper birational maps $f : T \to S$ one has $f_*[\O_T] =
[\O_S]$.  This is true if $T$ and $S$ have rational singularities,
which holds in all cases considered due to Lakshmibai and Magyar's
result \cite{lakshmibai.magyar:degeneracy}.


As mentioned above, the coefficients $c_\mu(r)$ depend only on the
side lengths of the rectangles $R_{ij}$, not on the integers $r_{ij}$
themselves.  Given that the coefficients have this property, they are
in fact uniquely given by the statement of \refthm{thm:structure} (see
\cite[\S2.2]{buch.fulton:chern}).  Regarding the signs of the
coefficients, we pose:

\begin{conj}
\label{conj:altsigns}
  The signs of the coefficients $c_\mu(r)$ alternate with the weight
  of $|\mu|$, {i.e.\ }$(-1)^{|\mu| - d(r)}\, c_\mu(r) \geq 0$.
\end{conj}

One particular case where this conjecture can be verified is when the
rectangle diagram only has two non-empty rows, {i.e.\ }$R_{ij}$ is empty
when $j-i > 2$.  This case includes all varieties of complexes.  When
all rectangles below the second row are empty, the inductive element
is given by $P_{\Bar r} = G_{R_{02}} \otimes G_{R_{13}} \otimes \dots
\otimes G_{R_{n-2,n}}$.

Now for a rectangular partition $R$, the coproduct constants
$d^R_{\sigma \tau}$ are given by the following simple rule.  Define a
{\em rook strip\/} to be a skew diagram which has at most one box in
any row or column.  Also, if $\tau$ is a partition which can be
contained in $R$, let $\Hat \tau$ denote $\tau$ rotated 180 degrees
and placed in the bottom-right corner of $R$.  We then have
\[ d^R_{\sigma \tau} = 
   \begin{cases} 
   (-1)^{|\sigma|+|\tau|-|R|} & \text{if $\sigma \cup \Hat \tau = R$
     and $\sigma \cap \Hat \tau$ is a rook strip;} \\
   0 & \text{otherwise.}
   \end{cases}
\]

\begin{cor}
  If $R_{ij}$ is empty for $j-i > 2$ then the formula is given by $P_r
  = \sum (-1)^{|\mu|-d(r)}\, G_{\mu_1} \otimes \dots \otimes
  G_{\mu_n}$, where the sum is over all partitions $\mu_i = (R_{i-1,i}
  + \sigma_i, \tau_{i-1})$, such that $\sigma_i \cup \Hat \tau_i =
  R_{i-1,i+1}$ and $\sigma_i \cap \Hat \tau_i$ is a rook strip for all
  $i$.
\end{cor}

Notice that the relations among the side lengths of the rectangles imply
that if $d^{R_{i-1,i+1}}_{\sigma_i, \tau_i} \neq 0$ then $\sigma_i$
always fits on the right side of $R_{i-1,i}$ and $\tau_i$ fits below
$R_{i,i+1}$, so the sequences of integers produced by the algorithm
are always partitions.

\[ \pic{40}{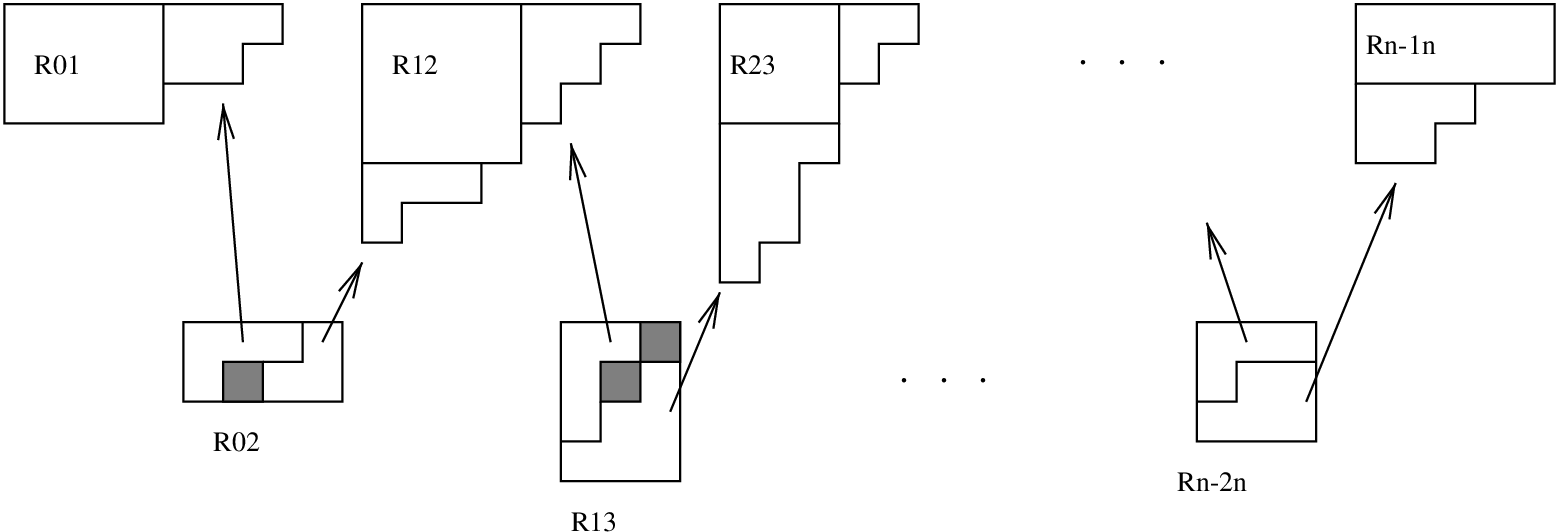} \]

In \cite{buch.fulton:chern} it is conjectured that the coefficients
$c_\mu(r)$ appearing in the cohomology formula (with $|\mu| = d(r)$)
are given as the number of sequences of semistandard Young tableaux
satisfying certain properties.  It would be very interesting to
generalize this conjecture to also give an expression for the more
general coefficients defined in this paper.


%% file: appgroth.tex
\section{Applications to Grothendieck polynomials}
\label{sec:appgroth}

In this section we will sketch how to apply our formula to give new
formulas for Grothendieck polynomials.  Our development is analogous
to \cite[\S 2.3]{buch.fulton:chern} and \cite{buch:stanley}.

Let $E_\bull$ be the sequence $F_1 \subset \dots \subset F_n \to H_n
\twoheadrightarrow \dots \twoheadrightarrow H_1$ considered in
Section~\ref{sec:groth}, and let $w \in S_{n+1}$ be a permutation.
Then $\Omega_w = \Omega_r(E_\bull)$ where $r = (r_{ij})$ are the
obvious rank conditions.
Set $x_i = 1 - L_i^{-1}$ and $y_i = 1 - M_i$ where $L_i = \ker(H_i \to
H_{i-1})$ and $M_i = F_i/F_{i-1}$.  The double Grothendieck polynomial
$\Groth_w(x;y)$ then becomes a special case of the quiver formula:
\[ \begin{split}
\Groth_w(x;y) &= [\O_{\Omega_r(E_\bull)}] \\
&= \sum c_\mu(r)\, G_{\mu_1}(F_2 - F_1) \cdots
   G_{\mu_{n-1}}(F_n - F_{n-1}) \cdot G_{\mu_n}(H_n - F_n) \cdot \\
& \hspace{2cm} G_{\mu_{n+1}}(H_{n-1} - H_n) \cdots G_{\mu_{2n-1}}(H_1 - H_2)
\end{split} \]
Set $\Tilde x_i = 1 - L_i = 1 - (1 - x_i)^{-1} = - (\sum_{k \geq 1}
x_i^k)$ and $\Tilde y_i = 1 - M_i^{-1} = - (\sum_{k \geq 1} y_i^k)$.
The formula can then be simplified using the identities
\[ G_\lambda(F_i - F_{i-1}) = G_\lambda(M_i) = 
   \begin{cases}
   (\Tilde y_i)^a & \text{if $\lambda = (a)$ is a row with $a$ boxes} \\
   0 & \text{otherwise}
   \end{cases}
\]
and
\[ G_\lambda(H_{i-1} - H_i) = G_\lambda(- L_i) =
   \begin{cases}
   (\Tilde x_i)^b & \text{if $\lambda = (1^b)$ is a column with $b$
   boxes} \\
   0 & \text{otherwise.}
   \end{cases}
\]
Using this we obtain a formula
\begin{equation}
\lab{eqn:newgroth}
\Groth_w(x;y) = \sum c_w(a,b,\lambda) \, \Tilde y_2^{a_2} \cdots
\Tilde y_n^{a_n} \, \Tilde x_2^{b_2} \cdots \Tilde x_n^{b_n} \,
G_\lambda(x; y) \,. 
\end{equation}
The sum is over exponents $a_2,\dots,a_n$ and $b_2,\dots,b_n$, and a
single partition $\lambda$, and $c_w(a,b,\lambda)$ is the coefficient
$c_\mu(r)$ for the sequence of partitions 
\[ \mu = ((a_2), \dots, (a_n), \lambda, 
   (1^{b_n}), \dots, (1^{b_2})) \,.
\]
Notice that it is not clear from the expression (\ref{eqn:newgroth})
that $\Groth_w(x;y)$ is a polynomial as opposed to a power series.

Now using the same arguments as in \cite{buch:stanley} we obtain
\[  \Groth_{1^m \times w}(x;y) = 
   \sum c_w(a,b,\lambda) \, \Tilde y_{2+m}^{a_2} \cdots \Tilde y_{n+m}^{a_n}
   \, \Tilde x_{2+m}^{b_2} \cdots \Tilde x_{n+m}^{b_n} \, G_\lambda(x;y)
\]
where $G_\lambda(x;y)$ is in variables $x_1,\dots,x_{n+m}$ and
$y_1,\dots,y_{n+m}$.  Letting $m$ tend to infinity in this expression,
it follows that
\[ G_w(x;y) = \sum_\lambda c_w(0,0,\lambda) \, G_\lambda(x;y) \,. \]

Thus we see that when the stable Grothendieck polynomial $G_w$ is
expressed in the basis $\{ G_\lambda \}$, the obtained coefficients
are special cases of the quiver coefficients $c_\mu(r)$ defined in
this paper.  Lascoux has recently shown that $(-1)^{|\lambda| -
  \ell(w)}\, c_w(0,0,\lambda) \geq 0$ \cite{lascoux:transition} which
confirms a special case of \refconj{conj:altsigns}.  In addition this
identity shows that the structure constants of $\Gamma$ are special
cases of the quiver coefficients $c_\mu(r)$, cf.\
\cite{buch:littlewood-richardson}.


%% file: jacobi-trudi.tex
\section{A generalized Jacobi-Trudi formula}
\label{sec:jacobi_trudi}

Recall that when one expands the Jacobi-Trudi determinant for the
Schur function $s_\lambda$ after the first row, one gets $s_\lambda =
\sum_{q \geq 0} (-1)^q \, s_{\lambda_1 + q} \cdot s_{\mu/(1^q)}$ where
$\mu = (\lambda_2,\lambda_3,\dots)$.  In this section we will prove a
generalization of this result for stable Grothendieck polynomials.

To state the formula in sufficient generality we need the following
definition.  If $I$ is a sequence of integers and $\lambda$ a
partition, we write
\[ G_{I\sslash \lambda} = 
   \sum_{\nu,\mu} \delta_{I,\nu} d^\nu_{\lambda \mu} \, G_\mu
   \in \Gamma \,.
\]
With this notation we have $\Delta G_I = \sum_\lambda G_{I\sslash
  \lambda} \otimes G_\lambda$.  Notice that when $I = \nu$ is a
partition, the element $G_{\nu\sslash \lambda}$ depends on both $\nu$
and $\lambda$ and not just the skew diagram $\nu/\lambda$ between
them.  For example $G_{\lambda\sslash \lambda} = 1$ if and only if
$\lambda$ is the empty partition.

\begin{thm}[Jacobi-Trudi formula]
\label{thm:jacobi_trudi}
If $a \in \Z$ is an integer and $I$ is a sequence of integers, then
\[ G_{a,I} = G_a \cdot G_I + \sum_{q \geq 1,\, t \geq 0}
   (-1)^q \binom{q-1+t}{t} \, G_{a+q+t} \cdot G_{I\sslash (1^q)} \,.
\]
\end{thm}

For proving this theorem, the following notation will get rid of a lot
of special cases.  We let $\bbn{n}{m}$ be the usual binomial
coefficient, except that we set $\bbn{-1}{0} = 1$:
\[ \bbn{n}{m} = \begin{cases} \binom{n}{m} &\text{if $0 \leq m \leq n$,} \\
   1 &\text{if $n = -1$ and $m = 0$,} \\ 
   0 &\text{otherwise.}
   \end{cases}
\]
\refthm{thm:jacobi_trudi} then asserts that $G_{a,I} = \sum_{q,t \geq
  0} (-1)^q \bbn{q-1+t}{t} G_{a+q+t} \cdot G_{I \sslash (1^q)}$, and
we have $\bbn{n}{m} = \bbn{n-1}{m-1} + \bbn{n-1}{m}$ whenever $m \leq
n$.

Notice that since $G_{a,I} = \sum_\mu \delta_{I,\mu} G_{a,\mu}$ and
$G_{I\sslash (1^q)} = \sum_\mu \delta_{I,\mu} G_{\mu\sslash (1^q)}$,
it is enough to prove the theorem in case $I = \mu$ is a partition.
We will give a bijective proof of the theorem when $a \geq \mu_1$.
For this we need the following combinatorial objects.  Recall that a
skew diagram is called a {\em horizontal strip\/} if no two boxes are
in the same column, and a {\em vertical strip\/} if no two boxes are
in the same row.  If both are true then the diagram is a rook strip.

\begin{defn}
A colored and marked Young diagram (CMYD) relative to a partition
$\mu$ is a quadruple of partitions $D = (\lambda_0 \subset \lambda
\subset \nu_0 \subset \nu)$ such that
\begin{romenum}
\item $\lambda \subset \mu$.
\item $\mu/\lambda_0$ is a vertical strip.
\item $\nu/\lambda$ is a horizontal strip.
\item $\lambda/\lambda_0$ and $\nu/\nu_0$ are both rook strips.
\item $\nu/\nu_0$ has no box in the top non-empty row of
  $\nu/\lambda$.
\end{romenum}
\end{defn}

We will regard a CMYD $D = (\lambda_0, \lambda, \nu_0, \nu)$ as the
Young diagram for $\nu$ in which the boxes of $\lambda$ are colored
white and the boxes of $\nu/\lambda$ are gray; the boxes in
$\lambda/\lambda_0$ and in $\nu/\nu_0$ are furthermore marked.
The axioms (i)--(v) then say that all white boxes are contained in
$\mu$; the boxes in $\mu$ which are not white form a vertical strip;
the gray boxes form a horizontal strip; the marked white boxes form a
rook strip and the marked gray boxes form a rook strip; and finally
the northernmost gray boxes are unmarked.  Let
\begin{align*}
g(D) &= \text{\# unmarked gray boxes in $D$} ~= |\nu_0/\lambda| \,, \\
w(D) &= \text{\# unmarked white boxes in $D$} ~= |\lambda_0| \,, \\
u(D) &= \text{\# unmarked boxes in $D$} ~= g(D) + w(D) \,\text{, and} \\
m(D) &= \text{\# marked boxes in $D$} ~= |\lambda/\lambda_0| +
|\nu/\nu_0| \,.
\end{align*}

We will write $G_D = G_\nu$.  By the coproduct Pieri rule of
\cite{buch:littlewood-richardson} or \refthm{thm:bialg} we have
$G_{\mu \sslash (1^q)} = \sum (-1)^{m(D)} G_D$, the sum over all CMYDs
relative to $\mu$ such that $w(D) = |\mu| - q$ and $D$ has no gray
boxes.  Then using Lenart's Pieri rule
\cite[Thm.~3.2]{lenart:combinatorial} or equation (\ref{eqn:lrmult})
we obtain
\begin{equation}
\label{eqn:cmyd_prod}
  G_p \cdot G_{\mu \sslash (1^q)} = \sum_D (-1)^{m(D)} G_D
\end{equation}
where this sum is over all CMYDs $D$ relative to $\mu$ with $g(D)
= p$ and $w(D) = |\mu| - q$.

\newcommand{\erow}[1]{\raisebox{18pt}{\picD{#1}}\hspace{8pt}}
\newcommand{\torow}[1]{\raisebox{9pt}{\picD{#1}}\hspace{8pt}}
\newcommand{\trerow}[1]{\picD{#1}\hspace{8pt}}
\begin{example}
\label{exm:cmyd}
  If we take $\mu = (1,1)$, $p = 2$, and $q = 1$, we have the
  following 6 CMYDs:
\[
\erow{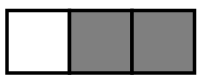} \torow{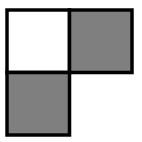} \torow{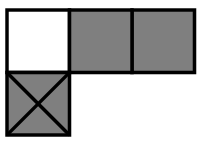} 
\torow{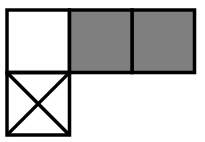} \trerow{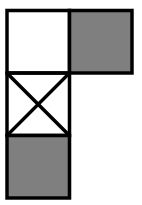} \trerow{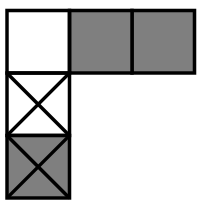} 
\]
It follows that $G_2 \cdot G_{1\,1 \sslash 1} = 
G_3 + G_{2\,1} - 2\,G_{3\,1} - G_{2\,1\,1} + G_{3\,1\,1}$.
\end{example}

Define the right vertical strip of $\mu$ to be the boxes in $\mu$ with
no boxes to the right of them.  We will say that a box in a CMYD $D$
is in $\mu$ resp.\ in the right vertical strip of $\mu$ if this is
true when the two diagrams are overlaid.  We will be interested in the
following four types of {\em special boxes\/} in $D$:

\begin{description}
\item[Type A] An unmarked gray box contained in $\mu$ which does not
  have a marked white box above it.

\item[Type B] Any white box (marked or unmarked) contained in the
  right vertical strip of $\mu$ which has no box under it.

\item[Type C] An unmarked gray box with a marked white box above it.
  
\item[Type D] A marked gray box such that the box above it is in the
  right vertical strip of $\mu$.
\end{description}

In \refexm{exm:cmyd} above, each diagram has exactly one special box.
From left to right, these boxes have types B, A, D, B, C, and D.  With
this notion we can rewrite the right hand side of the formula of
\refthm{thm:jacobi_trudi} as follows:

\begin{lemma}
\label{lemma:no_special}
For any partition $\mu$ and integer $a \geq \mu_1$ we have
\[ \sum_{q,t \geq 0} (-1)^q \bbn{q-1+t}{t} G_{a+q+t} \cdot G_{\mu\sslash (1^q)}
   = \sum_D (-1)^{|\mu| + w(D) + m(D)}
   \bbn{g(D)-a-1}{u(D)-a-|\mu|} G_D 
\]
where the sum is over all CMYDs relative to $\mu$ with no special
boxes.
\end{lemma}

\begin{proof}
  It follows from equation (\ref{eqn:cmyd_prod}) that the asserted
  identity is true if we sum over all CMYDs relative to $\mu$.  We
  will prove that the terms for which $D$ has special boxes cancel
  each other out in the right hand side.  Notice that each column of a
  CMYD can have at most one special box.  We will group each CMYD $D$
  for which the leftmost special box is of type A with two other CMYDs
  whose leftmost special boxes are of type B, such that the
  contributions from these three diagrams cancel.  Similarly a diagram
  with a leftmost special box of type C will be grouped with two
  diagrams with leftmost special boxes of type D.
  
  Notice that if $D$ is a CMYD relative to $\mu$ such that $u(D) - a -
  |\mu| \geq 0$ and $D$ contains a special box, then $D$ has at least
  $a+1$ gray boxes, so the top row of $D$ contains an unmarked gray
  box which is outside $\mu$.  Notice also that since $\bbn{n}{m} =
  \bbn{n-1}{m-1} + \bbn{n-1}{m}$ whenever $m \leq n$, we have
\begin{equation}
\label{eqn:special_cancel}
\bbn{g(D)-a-1}{u(D)-a-|\mu|} = 
   \bbn{g(D)-a-2}{u(D)-a-|\mu|-1} +
   \bbn{g(D)-a-2}{u(D)-a-|\mu|}
\end{equation}
for any diagram $D$ such that $w(D) < |\mu|$.

Now let $D$ be a CMYD whose leftmost special box is of type A.  The
conditions for a type A box then make it possible to change this box
into a white box or a marked white box, while the diagram continues to
be a CMYD.  Here it is important that the top row of $D$ contains at
least one unmarked gray boxes outside $\mu$, since this ensures that
the modified diagram satisfies axiom (v).
\[ \raisebox{-8pt}{\picB{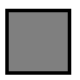}} ~~~~\longleftrightarrow~~~~ 
   \raisebox{-8pt}{\picB{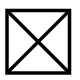}} ~~+~~ 
   \raisebox{-8pt}{\picB{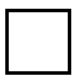}}
\]
The signs of the contributions from the two new diagrams are the
opposite of the sign of the contribution from $D$.  Since $w(D) <
|\mu|$, the contributions from all three diagrams therefore add to
zero by equation (\ref{eqn:special_cancel}).  Notice that any special
box of type B can be changed to a gray box.  Therefore all diagrams
with a leftmost special box of type B get canceled in this way.

Now suppose the leftmost special box in $D$ is of type C.  In this
case the box can be changed to a marked gray box while the box above
is either marked or unmarked white.
\[ \raisebox{-16pt}{\picB{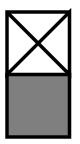}} ~~~~\longleftrightarrow~~~~ 
   \raisebox{-16pt}{\picB{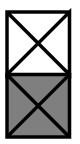}} ~~+~~ 
   \raisebox{-16pt}{\picB{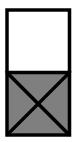}}
\]
Again the new diagrams give contributions of opposite sign from that
of $D$, and $w(D) < |\mu|$, so the contributions of all three diagrams
cancel by equation (\ref{eqn:special_cancel}).  Finally all diagrams
with a leftmost special box of type D are taken care of in this way,
since any diagram with a type D box can be changed so the special box
turns into type C.
\end{proof}

\begin{lemma}
\label{lemma:unique_cmyd}
Let $D$ be a CMYD with no special boxes and assume $u(D) \geq \mu_1 +
|\mu|$.  Then $D = (\mu, \mu, \nu, \nu)$ where $\nu = (g(D), \mu) =
(g(D), \mu_1, \mu_2, \dots)$.
\end{lemma}
\begin{proof}
  We start by observing that $D$ has no marked white boxes.  If $D$
  has such a box, then since it is not special, there must be a gray
  box below it.  But this gray box must then be special of type C or
  D, a contradiction.  Notice also that no unmarked gray boxes can be
  contained in $\mu$, since these would necessarily be special of type
  A.
  
  Now suppose $D$ contains a marked gray box, and consider the
  northernmost such box.  Since this box is not special (and not in
  the top row of $D$), the white box above it is not in the right
  vertical strip of $\mu$.  Now consider the row of boxes in $D$ to
  the right of this white box.  If this row contains a box in the
  right vertical strip of $\mu$, then this would necessarily be a
  special box of type B.  We conclude that if $D$ contains a marked
  gray box then some box northeast of this box is contained in $\mu$
  but not in $D$.
  
  Now assume that $\mu$ is not contained in $D$ and consider the
  northernmost row where $D$ is missing boxes from $\mu$.  Since $D$
  contains at least $\mu_1$ gray boxes, this can't be the top row, and
  the row above must contain a box in the right vertical strip of
  $\mu$ which has no box below it.  Since this box can't be marked
  gray by the argument above, it must be special of type $B$, again a
  contradiction.
  
  We conclude that $\mu$ is contained in $D$ and that all boxes from
  $\mu$ are white.  To prevent these white boxes from being special,
  there must furthermore be a gray box in each column of $D$.  This
  proves the result.
\end{proof}

The preceding two lemmas essentially prove \refthm{thm:jacobi_trudi}
when $I = \mu$ is a partition and $a \geq \mu_1$.  For the general
case of the theorem we will also need the following lemma.  Let
$h_i(x)$ denote the complete symmetric function of degree $i$.

\begin{lemma}
\label{lemma:gtos}
  For any integer $k \in \Z$ we have $G_k(x) = (1 - G_1(x)) \cdot
  \sum_{i \geq 0} h_{k+i}(x)$.
\end{lemma}
\begin{proof}
  If $k \geq 1$, it follows from \cite[Thm.~2.2]{lenart:combinatorial}
  that $G_k(x) = \sum_{p \geq 0} (-1)^p\, s_{(k,1^p)}(x)$.
  Alternatively this can be deduced from equation (\ref{eqn:single}),
  see e.g.\ \cite[\S6]{buch:littlewood-richardson}.  Notice in
  particular that $1 - G_1(x) = \sum_{p \geq 0} (-1)^p\, e_p(x)$.  For
  $k \geq 1$ the lemma therefore follows from the identity
  $\sum_{i=0}^p (-1)^i\, h_{k+i}\, e_{p-i} = s_{(k,1^p)}$.  When $k
  \leq 0$ the lemma is true because $\sum_{p \geq 0} (-1)^p\, e_p$ is
  the inverse power series to $\sum_{i \geq 0} h_i$.
\end{proof}

\begin{proof}[Proof of \refthm{thm:jacobi_trudi}]
  Suppose at first that $a \geq \mu_1$.  If $D$ is a CMYD relative to
  $\mu$ with no special boxes such that its coefficient
  $\bbn{g(D)-a-1}{u(D)-a-|\mu|}$ is non-zero, then since $u(D) \geq a
  + |\mu|$ we conclude by \reflemma{lemma:unique_cmyd} that $D =
  (\mu,\mu, \nu,\nu)$ where $\nu = (g(D),\mu)$.  But then we have
  $w(D) = |\mu|$ and $\bbn{g(D)-a-1}{g(D)-a} \neq 0$, so $g(D) = a$.
  The theorem therefore follows from \reflemma{lemma:no_special} in
  all cases where $a \geq \mu_1$.

For the general case it is enough to show that
\[ \GG_{a,\mu}(x_1,\dots,x_n) = \sum_{q,t \geq 0} (-1)^q
   \bbn{q-1+t}{t} \, G_{\mu\sslash (1^q)}(x_1,\dots,x_n) \cdot 
   G_{a+q+t}(x_1,\dots,x_n)
\]
where $n \geq 1 + \max(\ell(\mu), -a)$; this is sufficient because any
partition $\lambda$ such that $G_\lambda$ occurs in either side of the
claimed identity must have length at most $\ell(\mu) + 1$, and the
stable Grothendieck polynomials for partitions of such lengths are
linearly independent when applied to $n$ variables.  For the rest of
this proof we will let $x$ denote the $n$ variables $x_1,\dots,x_n$.

Let $\GG_\mu^{(i)}$ be the cofactor obtained by removing the first
row and the $i+1$'st column of the determinant defining
$\GG_{a,\mu}(x)$.  Notice that this does not depend on $a$, and we
have
\begin{equation}
\label{eqn:jt_lhs}
  \GG_{a,\mu}(x) = \sum_{i = 0}^{n-1} (-1)^i \,
  \GG_\mu^{(i)}(x) \cdot h_{a+i}(x) \,.
\end{equation}
Now using \reflemma{lemma:gtos} we obtain
\begin{equation}
\label{eqn:jt_rhs}
\begin{split}
 & \sum_{q,t \geq 0} (-1)^q \bbn{q-1+t}{t} \, 
   G_{\mu\sslash (1^q)}(x) \cdot G_{a+q+t}(x) \\
 &~~= \sum_{q,t \geq 0,\, i \geq q+t} (-1)^q \bbn{q-1+t}{t} \, 
   G_{\mu\sslash (1^q)}(x) \cdot (1-G_1(x)) \cdot h_{a+i}(x) \\
 &~~= \sum_{i \geq 0} \left( (1-G_1(x)) \cdot \sum_{q+t \leq i}
   (-1)^q \bbn{q-1+t}{t} G_{\mu\sslash (1^q)}(x) \right) \cdot h_{a+i}(x)
   \,.
\end{split} 
\end{equation}
Since (\ref{eqn:jt_lhs}) is equal to (\ref{eqn:jt_rhs}) for all $a
\geq \mu_1$, the theorem follows from the following lemma.
\end{proof}

\begin{lemma}
  Let $f_j \in \Z[[x_1,\dots,x_n]]$ be a power series for each $j \geq
  0$ and assume that
\begin{equation}
\label{eqn:hshift}
\sum_{j \geq 0} h_{a+j}(x_1,\dots,x_n) \cdot f_j = 0
\end{equation}
holds for all sufficiently large $a \in \N$.  Then (\ref{eqn:hshift})
is true for all $a \geq 1-n$.
\end{lemma}
\begin{proof}
  Since the form of each fixed degree in (\ref{eqn:hshift}) must be zero,
  we can assume that each $f_j$ is a polynomial and that $f_j = 0$ for $j
  > d$ for some $d \in \N$.  Assume at first that $d < n$ and let
  (\ref{eqn:hshift}) be true whenever $a \geq N$.  By assumption we
  then have
\[
\begin{bmatrix}
h_{N+d} & h_{N+d+1} & \dots & h_{N+2d} \\
h_{N+d-1} & h_{N+d} & \dots & h_{N+2d-1} \\
\vdots & \vdots && \vdots \\
h_N & h_{N+1} & \dots & h_{N+d}
\end{bmatrix} 
\begin{bmatrix}
f_0 \\ f_1 \\ \vdots \\ f_d
\end{bmatrix} =
\begin{bmatrix}
0 \\ 0 \\ \vdots \\ 0
\end{bmatrix} \,.
\]
Since the determinant of the matrix is the Schur polynomial
$s_{(N+d)^{d+1}}(x_1,\dots,x_n) \linebreak \neq 0$, we conclude that
each $f_j = 0$.

Now assume $d \geq n$.  If $a \geq 1-n$ then since $a + d \geq 1$  we
get
\[ h_{a+d}(x_1,\dots,x_n) = \sum_{j=0}^{d-1} (-1)^{d-j+1} 
h_{a+j}(x_1,\dots,x_n) e_{d-j}(x_1,\dots,x_n) \,. 
\]
So if we put $g_j = f_j + (-1)^{d-j+1} e_{d-j}(x_1,\dots,x_n) f_d$,
the left hand side of (\ref{eqn:hshift}) is equal to
\[ \sum_{j=0}^{d-1} h_{a+j}(x_1,\dots,x_n) \cdot g_j \,. \]
Since this is equal to zero for all large $a$, we conclude it is zero
for all $a \geq 1 - n$ by induction on $d$.
\end{proof}


%% file: gysin.tex
\section{A Gysin formula}
\label{sec:gysin}

In this section we prove of a $K$-theory parallel of a Gysin formula
of Pragacz \cite{pragacz:enumerative, fulton.pragacz:schubert}.  We
start with a Lemma which indicates that the classes $G_k(F)$ are the
right $K$-theoretic generalizations of {\em Segre classes\/} of a
vector bundle $F$.

\begin{lemma}
\label{lemma:segre}
Let $F$ be a vector bundle of rank $f$ over a variety $X$.  Let $\pi :
\P^*(F) \to X$ be the dual projective bundle of $F$ and let $Q$ be the
tautological quotient of $\pi^* F$.  Then for any $k \in \Z$ we have
$\pi_*(G_k(Q)) = G_{k-f+1}(F)$.
\end{lemma}
\begin{proof}
  This is clearly true if $k \leq 0$.  For $1 \leq k \leq f-1$ we
  write $G_k(Q) = (1-Q^{-1})^k = 1 + \sum_{i=1}^k (-1)^i Q^{-i}$, and
  the lemma follows because $R^j\pi_*(Q^{-i}) = 0$ for all $j$ and $1
  \leq i \leq f-1$.
  
  Finally, if $k \geq f$ we set $E = \O_X^{\oplus k}$ and form the
  bundle $H = \Hom(E,F) \to X$.  Then construct the fiber square:
\[ \xymatrix{
   Y \ar[r]^{\pi'} \ar[d] & H \ar[d] \\
   \P^*(F) \ar[r]^\pi & X} 
\]
We will suppress pullback notation for vector bundles.  By \cite[Lemma
1]{buch.fulton:chern} the locus $Z(E \to Q)$ in $Y$ is mapped
birationally onto $\Omega_{f-1}(E \to F) \subset H$.  Using
\refthm{thm:porteous} and the fact that determinantal varieties have
rational singularities \cite{gonciulea.lakshmibai:singular,
  lakshmibai.magyar:degeneracy} we therefore get
\[ \pi'_*(G_k(Q)) = \pi'_*( [\O_{Z(E \to Q)}] ) 
   = [ \O_{\Omega_{f-1}(E \to F)} ] 
   = G_{k-f+1}(F) \,.
\]
Since pullback along the vertical maps are isomorphisms which are
compatible with the horizontal pushforward maps, the lemma follows
from this.
\end{proof}

In the above lemma we have applied the usual proof of the
Thom-Porteous formula in the opposite way.  Alternatively one can
prove the identity $G_p(x_1,\dots,x_n) = 1 - \prod_{i=1}^n (1-x_i)
\cdot \sum_{j \geq 0} \binom{n+p-1}{n+j} h_j(x_1-1, \dots, x_n-1)$ by
induction on $n$.  The lemma follows from this using the identities
$\pi_*(Q^i) = S^i F$, $R^j \pi_*(Q^i) = 0$ for $0 < j < f-1$, and
$R^{f-1} \pi_*(Q^{-f-i}) = \det(F^\vee) \otimes S^i F^\vee$.

\begin{lemma}
\label{lemma:gysin}
Let $F$ and $E$ be vector bundles.  For any integers $m, i$ with $i
\geq \rank(E)$ we have $\sum_{k \geq 0} G_{m+k}(F) \, G_{i \sslash
  k}(-E) = G_{m+i}(F - E)$.
\end{lemma}
\begin{proof}
Notice that for any $p \in \Z$ we have $G_{p \sslash 0} = G_p$ and
$G_{p \sslash k} = G_{p-k} - G_{p-k+1}$ for all $k > 0$.  Since
$G_j(-E) = 0$ for $j > \rank(E)$ this means that $G_{p \sslash k}(-E)
= G_{p-k}(-E) - G_{p-k+1}(-E)$ if either $k > 0$, or $k \geq 0$ and $p
\geq \rank(E)$.

If $m \geq 0$ then $G_{m+i \sslash k}(-E) = 0$ for $k < m$ and $G_{m+i
  \sslash k}(-E) = G_{i \sslash k-m}(-E)$ for $k \geq m$, so
$G_{m+i}(F-E) = \sum_{k \geq 0} G_k(F)\, G_{m+i \sslash k}(-E) =
\sum_{k \geq 0} G_{m+k}(F)\, G_{i \sslash k}(-E)$ as required.

For $m \leq 0$ we have $\sum_{k = 0}^{-m} G_{m+k}(F)\, G_{i \sslash
  k}(-E) = \sum_{k=0}^{-m} G_{i \sslash k}(-E) = G_{m+i}(-E)$, so the
left hand side in the lemma is $G_{m+i}(-E) + \sum_{k \geq 1-m}
G_{m+k}(F)\, G_{i \sslash k}(-E) = G_{m+i}(-E) + \sum_{k \geq 1} G_k(F)\,
G_{m+i \sslash k}(-E) = G_{m+i}(F-E)$ as required.
\end{proof}

\begin{thm}
\label{thm:gysin}
Let $E$ and $F$ be bundles on $X$ of ranks $e$ and $f$.  Let $f = d +
q$ and let $\pi : \Gr(d, F) \to X$ be the Grassmann bundle of
$d$-planes in $F$ with universal exact sequence $0 \to A \to \pi^* F
\to Q \to 0$.  Let $I = (I_1, \dots, I_q)$ and $J = (J_1, J_2, \dots)$
be sequences of integers such that $I_j \geq e$ for all $j$.  Then
\[ \pi_*(G_I(Q - E) \cdot G_J(A - E)) = G_{I - (d)^q,J}(F - E) \,. \]
\end{thm}
\begin{proof}
Let $\Fl(d, f-1; F)$ be the variety of partial flags $A \subset H
\subset F$ such that $A$ has rank $d$ and $H$ has rank $f-1$.  Then
form the commutative diagram from \cite{jozefiak.lascoux.ea:classes}:
\[ \xymatrix{
   \Fl(d,f-1;F) \ar[r] \ar[d] & \Gr(d,F) \ar[d]^\pi \\
   \P^*(F) \ar[r] & X} 
\]
The formula can now be proved by calculating the pushforward to $X$ of
the class $G_{I_1+q-1}(F/H - E) \cdot G_{\Tilde I}(H/A - E) \cdot
G_J(A - E)$ in two different ways, using descending induction on $d$.
Here $\Tilde I$ is the sequence $(I_2, \dots, I_q)$.  

We are therefore reduced to the case $d = f - 1$ where $\Gr(d,F) =
\P^*(F)$ is a projective bundle.  Notice that since $Q$ is now a line
bundle we have $G_k(Q) = (1 - [Q^\vee])^k$ for $k \geq 0$.  Using this
we get the following identities in $K^\circ \P^*(F)$.

\[ \begin{split}
& G_i(Q-E) \cdot G_J(A-E) = G_i(Q-E) \cdot G_J(F - E \oplus Q) \\
&~~= \sum_{k \geq 0} G_k(Q) \cdot G_{i \sslash k}(-E) \cdot 
     \sum_{\ell \geq 0} G_{(1^\ell)}(-Q) \cdot G_{J\sslash(1^\ell)}(F-E) \\
&~~= \sum_{k,\ell \geq 0} G_k(Q) \cdot G_\ell(Q^\vee)
     \cdot G_{i\sslash k}(-E) \cdot G_{J\sslash (1^\ell)}(F-E) \\
&~~= \sum_{k,\ell \geq 0} (-1)^\ell [Q]^\ell G_{k+\ell}(Q) \cdot
      G_{i\sslash k}(-E) \cdot G_{J\sslash (1^\ell)}(F-E) \\
&~~= \sum_{k,\ell \geq 0} (-1)^\ell (1-G_1(Q))^{-\ell}
     \cdot G_{k+\ell}(Q) \cdot G_{i\sslash k}(-E) \cdot G_{J\sslash (1^\ell)}(F-E)\\
&~~= \sum_{k,\ell \geq 0} (-1)^\ell \sum_{t \geq 0}
     \bbn{\ell-1+t}{t} G_1(Q)^t \cdot G_{k+\ell}(Q) \cdot G_{i\sslash k}(-E) \cdot
     G_{J\sslash (1^\ell)}(F-E) \\
&~~= \sum_{k,\ell,t \geq 0} (-1)^\ell \bbn{\ell-1-t}{t}
     G_{k+\ell+t}(Q) \cdot G_{i\sslash k}(-E)\cdot G_{J\sslash (1^\ell)}(F-E) \\
\end{split} \]
The step replacing $(1 - G_1(Q))^{-\ell}$ with its power series
expansion is valid since $G_1(Q)^t$ is zero for $t > \dim \P^*(F)$.
Now using \reflemma{lemma:segre}, \reflemma{lemma:gysin}, and
\refthm{thm:jacobi_trudi} we get
\[ \begin{split}
& \pi_*( G_i(Q - E) \cdot G_J(A - E)) \\
&~~= \sum_{k,\ell,t \geq 0} (-1)^\ell \bbn{\ell-1+t}{t}
    G_{k+\ell+t-d}(F) \cdot G_{i\sslash k}(-E) \cdot 
    G_{J\sslash (1^\ell)}(F-E) \\
&~~= \sum_{\ell,t \geq 0} (-1)^\ell \bbn{\ell-1+t}{t}
     G_{\ell+t+i-d}(F-E)\cdot G_{J\sslash (1^\ell)}(F-E) \\
&~~= G_{i-d,J}(F-E)
\end{split} \]
which is what we want to prove.
\end{proof}

Continuing the remark after \reflemma{lemma:segre}, notice that
\refthm{thm:porteous} is a consequence of \refthm{thm:gysin} once we
prove that the structure sheaf of a zero section $Z(E \to F)$ is given
by $G_{(e)^f}(F - E)$, see e.g.\ \cite[\S 14.4]{fulton:intersection}.
This in turn follows from
\cite[eqn.~(7.1)]{buch:littlewood-richardson}.


We will finish this paper with the following somewhat surprising
consequence of \refthm{thm:gysin}.

\begin{cor}
  Let $\pi: \Gr(d,F) \to X$ be a Grassmann bundle with tautological
  quotient bundle $Q$ of rank $q$.  For partitions $\lambda$ and $\mu$
  of lengths at most $q$, we have
\[ \pi_*(G_\lambda(Q)) \cdot \pi_*(G_\mu(Q)) =
   \pi_*(G_\lambda(Q) \cdot G_\mu(Q)) + \sum_{\ell(\nu) > q}
   c^\nu_{\lambda \mu} \, G_{\Tilde \nu}(F)
\]
where the sum is over all partitions $\nu$ of length strictly greater
than $q$ and $\Tilde \nu$ denotes $\nu$ with the first $d$ columns
removed.  In particular, if $\ell(\lambda) + \ell(\mu) \leq q$ then we
get $\pi_*(G_\lambda(Q) \cdot G_\mu(Q)) = \pi_*(G_\lambda(Q)) \cdot
\pi_*(G_\mu(Q))$.
\end{cor}

In other words, $\pi_*$ behaves like a ring homomorphism for short
partitions!

\begin{proof}
  In \cite[\S 7]{buch:littlewood-richardson} it is shown that the
  linear map $\Gamma \to \Gamma$ defined by $G_\nu \mapsto G_{\Tilde
    \nu}$ is a ring homomorphism.  Using this we get
  $\pi_*(G_\lambda(Q)) \cdot \pi_*(G_\mu(Q)) = G_{\Tilde \lambda}(F)
  \cdot G_{\Tilde \mu}(F) = \sum_\nu c^\nu_{\lambda \mu} G_{\Tilde
    \nu}(F)$.  On the other hand we have $G_\lambda(Q) \cdot
  G_\mu(Q) = \sum_{\ell(\nu) \leq q} c^\nu_{\lambda \mu} G_\nu(Q)$,
  so $\pi_*(G_\lambda(Q) \cdot G_\mu(Q)) = \sum_{\ell(\nu) \leq q}
  c^\nu_{\lambda \mu} G_{\Tilde \nu}(F)$.  The corollary follows from
  this.
\end{proof}


%% file: bibliography.tex
\providecommand{\bysame}{\leavevmode\hbox to3em{\hrulefill}\thinspace}
